\documentclass[12pt, a4paper, reqno]{amsart}
\usepackage{amsmath}
\usepackage{amsfonts}
\usepackage{amssymb}
\usepackage{color}
\usepackage{comment,graphicx}
\usepackage[T1]{fontenc}
\usepackage{url}
\usepackage{ae}
\usepackage{mathtools}
\usepackage{esint}
\usepackage{latexsym}
\usepackage{amscd}

\def\phi{\varphi}
\def\rho{\varrho}
\def\epsilon{\varepsilon}
\numberwithin{equation}{section}
\theoremstyle{plain}
\newtheorem{theorem}[equation]{Theorem}
\newtheorem{lemma}[equation]{Lemma}

\theoremstyle{definition}
\newtheorem{definition}[equation]{Definition}

\theoremstyle{remark}
\newtheorem{remark}[equation]{Remark}
\newtheorem{example}[equation]{Example}

\renewcommand{\leq}{\leqslant}
\renewcommand{\geq}{\geqslant}

\input{tcilatex}

\begin{document}
\title[Dilation in function spaces with general weights]{Dilation in
function spaces with general weights}
\author[D. Drihem]{Douadi Drihem}
\address{Douadi Drihem\\
M'sila University\\
Department of Mathematics\\
Laboratory of Functional Analysis and Geometry of Spaces\\
M'sila 28000, Algeria.}
\email{douadidr@yahoo.fr, douadi.drihem@univ-msila.dz}
\thanks{ }
\date{\today }

\begin{abstract}
In this paper, we present more regularity conditions which ensure the
boundedness of dilation operators on \ Besov and Triebel-Lizorkin spaces
equiped with general weights.
\end{abstract}

\maketitle

\section{Introduction}

Function spaces play a central role in mathematical analysis especially in
partial differential equations. Some example of these spaces can be
mentioned such as: Besov and Triebel-Lizorkin spaces. The theory of these
spaces had a remarkable development in part due to its usefulness in
applications. We refer the reader to the monographs \cite{T1} and \cite{T2}
for further details, historical remarks and more references on these spaces.

In recent years many researchers have modified the classical spaces and have
generalized the classical results to these modified ones. For example:
function spaces of generalized smoothness. These type of function spaces
have been introduced by several authors. We refer, for instance, to Cobos
and Fernandez \cite{CF88}, Goldman \cite{Go79} and \cite{Go83}, and Kalyabin 
\cite{Ka83}; see also\ Besov \cite{B03} and \cite{B05}, and Kalyabin and
Lizorkin \cite{Kl87}.

These\ type of function\ spaces appear in the study of trace spaces on
fractals, see Edmunds and Triebel \cite{ET96} and \cite{ET99}, were they
introduced the spaces $B_{p,q}^{s,\Psi }$, where $\Psi $ is a so-called
admissible function, typically of log-type near $0$. For a complete
treatment of these spaces we refer the reader the work of Moura \cite{Mo01}.
More general function spaces of generalized smoothness can be found in
Farkas and Leopold \cite{FL06}, and reference therein.

A. Tyulenev has been introduced in \cite{Ty15} a new family of Besov space
of variables smoothness which cover many classes of Besov spaces. Based on
this new weighted class and the Littlewood-Paley theory the author in \cite%
{D6} and \cite{D6.1} has been introduced the function spaces $B_{p,q}(%
\mathbb{R}^{n},\{t_{k}\}_{k\in \mathbb{N}_{0}})$ and $F_{p,q}(\mathbb{R}%
^{n},\{t_{k}\}_{k\in \mathbb{N}_{0}})$ which cover weighted Besov and
Triebel-Lizorkin spaces, respectively. Several results, concerning, for
instance, Sobolev embeddings, atomic, molecular and wavelet decompositions
are presented.

The purpose of the present paper is to study the dilation operators in $%
A_{p,q}(\mathbb{R}^{n},\{t_{k}\}_{k\in \mathbb{N}_{0}})$ spaces, where\ we
use this notation to denote either $B_{p,q}(\mathbb{R}^{n},\{t_{k}\}_{k\in 
\mathbb{N}_{0}})$ or $F_{p,q}(\mathbb{R}^{n},\{t_{k}\}_{k\in \mathbb{N}%
_{0}}) $. Their behavior is well known in Besov and Triebel-Lizorkin spaces,
see \cite[3.4.1]{T1}. Further results can be found in \cite{Sc09}, \cite%
{ScV09} and \cite{V08}. Allowing $\{t_{k}\}_{k\in \mathbb{N}_{0}}$ to vary
from point to point will raise extra difficulties which, in general, are
overcome by imposing some regularity assumptions on this smoothness. By
these additional assumptions we ensure the boundedness of these operators on 
$A_{p,q}(\mathbb{R}^{n},\{t_{k}\}_{k\in \mathbb{N}_{0}})$ spaces, but with
some appropriate assumptions. More precisely we shall show the following
result:

\begin{theorem}
\label{Dilation}\textit{Let }$1\leq p<\infty ,1\leq q<\infty ,\alpha
_{1},\alpha _{2}\in \mathbb{R},\alpha =(\alpha _{1},\alpha _{2}),\lambda
\geq 1$. Let $\alpha _{2}\geq \alpha _{1}>0$. There exists $1\leq \theta
<p<\infty \ $such that for all $\{t_{k}\}_{k\in \mathbb{N}_{0}}\in X_{\alpha
,\sigma ,p}$, with $\sigma =(\sigma _{1}=\theta \left( \frac{p}{\theta }%
\right) ^{\prime },\sigma _{2}\geq p)$ and all $f\in A_{p,q}(\mathbb{R}%
^{n},\{t_{k}\}_{k\in \mathbb{N}_{0}})$%
\begin{equation}
\big\|f(\lambda \cdot )\big\|_{A_{p,q}(\mathbb{R}^{n},\{t_{k}\}_{k\in 
\mathbb{N}_{0}})}\leq c\lambda ^{\alpha _{2}-\frac{n}{p}}H\big\|f\big\|%
_{A_{p,q}(\mathbb{R}^{n},\{t_{k}\}_{k\in \mathbb{N}_{0}})},  \label{key-est}
\end{equation}%
where%
\begin{equation*}
H=\sup_{k\geq i,m\in \mathbb{Z}^{n}}\frac{\big\|t_{k-i}(\lambda ^{-1}\cdot )%
\big\|_{L_{p}(Q_{k-i,m})}}{\big\|t_{k-i}\big\|_{L_{p}(Q_{k-i,m})}},
\end{equation*}%
the constant $c>0$ independent of $\lambda $ and $\lambda <2^{i}\leq
2\lambda .$
\end{theorem}

As a consequence, our result cover the classical case, see \cite[3.4.1]{T1},
also for Besov and Triebel-Lizorkin spaces equipped with power weights.
Concerning Sobolev spaces $W_{p}^{k}(\mathbb{R}^{n},w)$ it holds%
\begin{equation*}
W_{p}^{k}(\mathbb{R}^{n},w)=F_{p}^{k}(\mathbb{R}^{n},w),\quad 1<p<\infty
,k\in \mathbb{N}_{0},w\in A_{p}(\mathbb{R}^{n}),
\end{equation*}%
see \cite[Theorem 2.8]{Bui82}, where $A_{p}(\mathbb{R}^{n})$ are the
Muckenhoupt classes, see Section 2. We can easily prove that%
\begin{equation}
\big\|f(\lambda \cdot )\big\|_{W_{p}^{k}(\mathbb{R}^{n},w)}\leq \lambda ^{k-%
\frac{n}{p}}\sup_{x\in \mathbb{R}^{n}}\frac{\omega (\lambda ^{-1}x)}{\omega
(x)}\big\|f\big\|_{W_{p}^{k}(\mathbb{R}^{n},w)},  \label{key-est1}
\end{equation}%
for all $f\in W_{p}^{k}(\mathbb{R}^{n},w)$. In Section 4 we prove that our
estimate $\mathrm{\eqref{key-est}}$ is better than $\mathrm{\eqref{key-est1}.%
}$

We mention that the boundedness of these operators in function spaces play
an important role in mathematical analysis. They appear in the
Gagliardo-Nirenberg inequalities \cite[Chapter 4]{T3} in the boundedness
properties of pseudodifferential operators on Besov spaces and
Triebel-Lizorkin \cite{Naibo}, and in the spline representations of Besov
and Triebel-Lizorkin spaces \cite{Si}.

\section{Basic tools}

Throughout this paper, we denote by $\mathbb{R}^{n}$ the $n$-dimensional
real Euclidean space, $\mathbb{N}$ the collection of all natural numbers and 
$\mathbb{N}_{0}=\mathbb{N}\cup \{0\}$. The letter $\mathbb{Z}$ stands for
the set of all integer numbers.\ The expression $f\lesssim g$ means that $%
f\leq c\,g$ for some independent constant $c$ (and non-negative functions $f$
and $g$).\ As usual for any $x\in \mathbb{R}$, $\left\lfloor x\right\rfloor $
stands for the largest integer smaller than or equal to $x$.\vskip5pt

By supp $f$ we denote the support of the function $f$, i.e., the closure of
its non-zero set. If $E\subset {\mathbb{R}^{n}}$ is a measurable set, then $%
|E|$ stands for the (Lebesgue) measure of $E$ and $\chi _{E}$ denotes its
characteristic function. By $c$ we denote generic positive constants, which
may have different values at different occurrences. \vskip5pt

A weight is a nonnegative locally integrable function on $\mathbb{R}^{n}$
that takes values in $(0,\infty )$ almost everywhere. For measurable set $%
E\subset \mathbb{R}^{n}$ and a weight $\gamma $, $\gamma (E)$ denotes 
\begin{equation*}
\int_{E}\gamma (x)dx.
\end{equation*}%
Given a measurable set $E\subset \mathbb{R}^{n}$ and $0<p\leq \infty $, we
denote by $L_{p}(E)$ the space of all functions $f:E\rightarrow \mathbb{C}$
equipped with the finite quasi-norm 
\begin{equation*}
\left\Vert f\right\Vert _{L_{p}(E)}:=\Big(\int_{E}\left\vert f(x)\right\vert
^{p}dx\Big)^{1/p}<\infty ,\text{\quad }0<p<\infty ,
\end{equation*}%
\begin{equation*}
\left\Vert f\right\Vert _{L_{\infty }(E)}:=\underset{x\in E}{\text{ess-sup}}%
\left\vert f(x)\right\vert <\infty .
\end{equation*}%
For a function $f$ in $L_{1}^{\mathrm{loc}}$, we set 
\begin{equation*}
M_{A}(f):=\frac{1}{|A|}\int_{A}\left\vert f(x)\right\vert dx
\end{equation*}%
for any $A\subset \mathbb{R}^{n}$. Furthermore, we put%
\begin{equation*}
M_{A,p}(f):=\Big(\frac{1}{|A|}\int_{A}\left\vert f(x)\right\vert ^{p}dx\Big)%
^{\frac{1}{p}},\quad 0<p<\infty .
\end{equation*}%
Further, given a measurable set $E\subset \mathbb{R}^{n}$ and a weight $%
\gamma $, we denote the space of all functions $f:\mathbb{R}^{n}\rightarrow 
\mathbb{C}$ with finite quasi-norm 
\begin{equation*}
\left\Vert f\right\Vert _{L_{p}(\mathbb{R}^{n},\gamma )}=\left\Vert f\gamma
|\right\Vert _{L_{p}(\mathbb{R}^{n})}
\end{equation*}%
by $L_{p}(\mathbb{R}^{n},\gamma )$.

Let $0<p\leq \infty $ and $0<q\leq \infty $. The space $L_{p}(\ell _{q})$ is
defined to be the set of all sequences $\{f_{k}\}_{k\in \mathbb{Z}}$ of
functions such that%
\begin{equation*}
\big\|\{f_{k}\}_{k\in \mathbb{N}_{0}}\big\|_{L_{p}(\ell _{q})}:=\big\|%
\left\Vert \{f_{k}\}_{k\in \mathbb{Z}}\right\Vert _{\ell _{q}}\big\|_{L_{p}(%
\mathbb{R}^{n})}<\infty .
\end{equation*}%
In the limiting case $q=\infty $ the usual modification is required. If $%
1\leq p\leq \infty $ and $\frac{1}{p}+\frac{1}{p^{\prime }}=1$, then $%
p^{\prime }$ is called the conjugate exponent of $p$.

The symbol $\mathcal{S}(\mathbb{R}^{n})$ is used in place of the set of all
Schwartz functions on $\mathbb{R}^{n}$. We define the Fourier transform of a
function $f\in \mathcal{S}(\mathbb{R}^{n})$ by 
\begin{equation*}
\mathcal{F}(f)(\xi ):=(2\pi )^{-n/2}\int_{\mathbb{R}^{n}}e^{-ix\cdot \xi
}f(x)dx,\quad \xi \in \mathbb{R}^{n}.
\end{equation*}%
In what follows, $Q$ will denote an cube in the space $\mathbb{R}^{n}$\ with
sides parallel to the coordinate axes and $l(Q)$\ will denote the side
length of the cube $Q$. For $k\in \mathbb{N}_{0}$ and $m\in \mathbb{Z}^{n}$,
denote by $Q_{k,m}$ the dyadic cube,%
\begin{equation*}
Q_{k,m}:=2^{-k}([0,1)^{n}+m).
\end{equation*}%
For the collection of all such cubes we use $\mathcal{Q}:=\{Q_{k,m}:k\in 
\mathbb{N}_{0},m\in \mathbb{Z}^{n}\}$. For each cube $Q$, we denote by $%
x_{k,m}$ the lower left-corner $2^{-k}m$ of $Q=Q_{k,m}$.

\subsection{Muckenhoupt weights}

The purpose of this subsection is to review some known properties of\
Muckenhoupt classes.

\begin{definition}
Let $1<p<\infty $. We say that a weight $\gamma $ belongs to the Muckenhoupt
class $A_{p}(\mathbb{R}^{n})$ if there exists a constant $C>0$ such that for
every cube $Q$ the following inequality holds 
\begin{equation}
M_{Q}(\gamma )M_{Q,\frac{p^{\prime }}{p}}(\gamma ^{-1})\leq C.
\label{Ap-constant}
\end{equation}
\end{definition}

The smallest constant $C$ for which $\mathrm{\eqref{Ap-constant}}$ holds,
denoted by $A_{p}(\gamma )$. As an example, we can take%
\begin{equation*}
\gamma (x)=|x|^{\alpha },\quad \alpha \in \mathbb{R}.
\end{equation*}%
Then $\gamma \in A_{p}(\mathbb{R}^{n})$, $1<p<\infty $, if and only if $%
-n<\alpha <n(p-1)$.

For $p=1$ we rewrite the above definition in the following way.

\begin{definition}
We say that a weight $\gamma $ belongs to the Muckenhoupt class $A_{1}(%
\mathbb{R}^{n})$ if there exists a constant $C>0$ such that for every cube $%
Q $\ and for a.e.\ $y\in Q$ the following inequality holds 
\begin{equation}
M_{Q}(\gamma )\leq C\gamma (y).  \label{A1-constant}
\end{equation}
\end{definition}

The smallest constant $C$ for which $\mathrm{\eqref{A1-constant}}$ holds,
denoted by $A_{1}(\gamma )$. The above classes have been first studied by
Muckenhoupt\ \cite{Mu72} and use to characterize the boundedness of the
Hardy-Littlewood maximal function on $L^{p}(\gamma )$, see the monographs 
\cite{GR85} and \cite{L. Graf08}\ for a complete account on the theory of
Muckenhoupt weights.

We recall a few basic properties of the class $A_{p}(\mathbb{R}^{n})$
weights, see \cite{L. Graf08}.

\begin{lemma}
\label{Ap-Property} Let $1\leq p<\infty $.\newline
$\mathrm{(i)}$ If $\gamma \in A_{p}(\mathbb{R}^{n})$, then for any $1\leq
p<q $, $\gamma \in A_{q}(\mathbb{R}^{n})$.\newline
$\mathrm{(ii)}$ Let $1<p<\infty $. $\gamma \in A_{p}(\mathbb{R}^{n})$ if and
only if $\gamma ^{1-p^{\prime }}\in A_{p^{\prime }}(\mathbb{R}^{n})$.\newline
$\mathrm{(iii)}$ Let $1\leq p<\infty $ and $\gamma \in A_{p}(\mathbb{R}^{n})$%
. There is $C>0$ such that for any cube $Q$ and a measurable subset $%
E\subset Q$%
\begin{equation*}
\left( \frac{|E|}{|Q|}\right) ^{p-1}M_{Q}(\gamma )\leq CM_{E}(\gamma ).
\end{equation*}%
$\mathrm{(iv)}$ Suppose that $\gamma \in A_{p}(\mathbb{R}^{n})$ for some $%
1<p<\infty $. Then there exists a $1<p_{1}<p<\infty $ such that $\gamma \in
A_{p_{1}}(\mathbb{R}^{n})$.\newline
$\mathrm{(v)}$ Let $\gamma \in A_{p}(\mathbb{R}^{n})$. Then $\gamma (\lambda
\cdot )\in A_{p}(\mathbb{R}^{n})$ for any $\lambda >0.$
\end{lemma}

\subsection{The weight class $X_{\protect\alpha ,\protect\sigma ,p}$}

If given $0<p\leq \infty $ a sequence of weights $\{t_{k}\}_{k\in \mathbb{N}%
_{0}}$ is such that $t_{k}\in L_{p}^{\mathrm{loc}}$ for $k\in \mathbb{N}_{0}$%
, then the weight sequence $\{t_{k}\}_{k\in \mathbb{N}_{0}}$ will be called
a $p$-admissible weight sequence. For a $p$-admissible weight sequence $%
\{t_{k}\}_{k\in \mathbb{N}_{0}}$\ we set%
\begin{equation*}
t_{k,m}:=\big\|t_{k}\big\|_{L_{p}(Q_{k,m})},\quad k\in \mathbb{N}_{0},m\in 
\mathbb{Z}^{n}.
\end{equation*}

Tyulenev\ \cite{Ty14} introduce the following new weighted class\ and use to
study Besov spaces of variable smoothness.

\begin{definition}
\label{Tyulenev-class}Let $0<p\leq \infty ,\alpha _{1}$, $\alpha _{2}\in 
\mathbb{R}$, $\sigma _{1}$, $\sigma _{2}$ $\in (0,+\infty ]$, $\alpha
=(\alpha _{1},\alpha _{2})$ and let $\sigma =(\sigma _{1},\sigma _{2})$. We
let $X_{\alpha ,\sigma ,p}=X_{\alpha ,\sigma ,p}(\mathbb{R}^{n})$ denote the
set of $p$-admissible weight sequences $\{t_{k}\}_{k\in \mathbb{N}_{0}}$
satisfying the following conditions. There exist numbers $C_{1},C_{2}>0$
such that for any $0\leq k\leq j$\ and every cube $Q,$%
\begin{eqnarray}
M_{Q,p}(t_{k})M_{Q,\sigma _{1}}(t_{j}^{-1}) &\leq &C_{1}2^{\alpha _{1}(k-j)},
\label{Asum1} \\
M_{Q,p}^{-1}(t_{k})M_{Q,\sigma _{2}}(t_{j}) &\leq &C_{2}2^{\alpha _{2}(j-k)}.
\label{Asum2}
\end{eqnarray}
\end{definition}

The constants $C_{1},C_{2}>0$ are independent of both the indexes $k$ and $j$%
.

\begin{remark}
We would like to mention that if $\{t_{k}\}_{k\in \mathbb{N}_{0}}$
satisfying $\mathrm{\eqref{Asum1}}$ with $\sigma _{1}=r\left( \frac{p}{r}%
\right) ^{\prime }$ and $0<r<p\leq \infty $, then $t_{k}^{p}\in A_{\frac{p}{r%
}}(\mathbb{R}^{n})$ for any $k\in \mathbb{N}_{0}$ with\ $0<r<p<\infty $ and $%
t_{k}^{-r}\in A_{1}(\mathbb{R}^{n})$ for any $k\in \mathbb{N}_{0}$\ with\ $%
p=\infty $.\newline
We say that $t_{k}\in A_{p}(\mathbb{R}^{n})$,\ $k\in \mathbb{N}_{0}$, $%
1<p<\infty $ have the same Muckenhoupt constant if%
\begin{equation*}
A_{p}(t_{k})=c,\quad k\in \mathbb{N}_{0},
\end{equation*}%
where $c$ is independent of $k$.
\end{remark}

\begin{example}
\label{Example1}Let $0<r<p<\infty $, a weight $\omega ^{p}\in A_{\frac{p}{r}%
}(\mathbb{R}^{n})$ and $\{s_{k}\}_{k\in \mathbb{N}_{0}}=\{2^{ks}\omega
^{p}(2^{-k})\}_{k\in \mathbb{N}_{0}}$, $s\in \mathbb{R}$. Clearly, $%
\{s_{k}\}_{k\in \mathbb{N}_{0}}$ lies in $X_{\alpha ,\sigma ,p}$ for $\alpha
_{1}=\alpha _{2}=s$, $\sigma =(r(\frac{p}{r})^{\prime },p)$.\ An example
illustrating the advantage of Definition \ref{Tyulenev-class} is given in 
\cite{Ty15}, \cite{Ty-N-L} and \cite{Ty-151}.
\end{example}

\begin{remark}
\label{Tyulenev-class-properties}Let $0<\theta \leq p\leq \infty $. Let $%
\alpha _{1}$, $\alpha _{2}\in \mathbb{R}$, $\sigma _{1},\sigma _{2}\in
(0,+\infty ]$, $\sigma _{2}\geq p$, $\alpha =(\alpha _{1},\alpha _{2})$ and
let $\sigma =(\sigma _{1}=\theta \left( \frac{p}{\theta }\right) ^{\prime
},\sigma _{2})$. Let a $p$-admissible weight sequence $\{t_{k}\}_{k\in 
\mathbb{N}_{0}}\in X_{\alpha ,\sigma ,p}$. Then%
\begin{equation*}
\alpha _{2}\geq \alpha _{1}.
\end{equation*}
\end{remark}

No we recall the vector-valued maximal inequality of Fefferman and Stein 
\cite{FeSt71}. As usual, we put%
\begin{equation*}
\mathcal{M(}f)(x):=\sup_{Q}\frac{1}{|Q|}\int_{Q}\left\vert f(y)\right\vert
dy,\quad f\in L_{1}^{\mathrm{loc}},
\end{equation*}%
where the supremum\ is taken over all cubes with sides parallel to the axis
and $x\in Q$. Also we set 
\begin{equation*}
\mathcal{M}_{\sigma }(f):=\left( \mathcal{M(}\left\vert f\right\vert
^{\sigma })\right) ^{\frac{1}{\sigma }},\quad 0<\sigma <\infty .
\end{equation*}

\begin{theorem}
Let $0<p<\infty ,0<q\leq \infty $ and $0<\sigma <\min (1,p,q)$. Then%
\begin{equation}
\Big\|\Big(\sum\limits_{k=0}^{\infty }\big(\mathcal{M}_{\sigma }(f_{k})\big)%
^{q}\Big)^{\frac{1}{q}}\Big\|_{L_{p}(\mathbb{R}^{n})}\lesssim \Big\|\Big(%
\sum\limits_{k=0}^{\infty }\left\vert f_{k}\right\vert ^{q}\Big)^{\frac{1}{q}%
}\Big\|_{L_{p}(\mathbb{R}^{n})}  \label{Fe-St71}
\end{equation}%
holds for all sequence of functions $\{f_{k}\}_{k\in \mathbb{N}_{0}}\in
L_{p}(\ell _{q})$.
\end{theorem}

We will make use of the following statement, see \cite{D6.1}.

\begin{lemma}
\label{key-estimate1}Let $1<\theta \leq p<\infty $\ and $1<q<\infty $. Let $%
\{t_{k}\}_{k\in \mathbb{N}_{0}}$\ be a $p$-admissible\ weight\ sequence\
such that $t_{k}^{p}\in A_{\frac{p}{\theta }}(\mathbb{R}^{n})$, $k\in 
\mathbb{N}_{0}$. Assume that $t_{k}^{p}$,\ $k\in \mathbb{N}_{0}$ have the
same Muckenhoupt constant, $A_{\frac{p}{\theta }}(t_{k})=c,k\in \mathbb{N}%
_{0}$. Then%
\begin{equation*}
\Big\|\Big(\sum\limits_{k=0}^{\infty }t_{k}^{q}\big(\mathcal{M(}f_{k})\big)%
^{q}\Big)^{\frac{1}{q}}\Big\|_{L_{p}(\mathbb{R}^{n})}\lesssim \Big\|\Big(%
\sum\limits_{k=0}^{\infty }t_{k}^{q}\left\vert f_{k}\right\vert ^{q}\Big)^{%
\frac{1}{q}}\Big\|_{L_{p}(\mathbb{R}^{n})}
\end{equation*}%
holds for all sequence of functions $\{f_{k}\}_{k\in \mathbb{N}_{0}}\in
L_{p}(\ell _{q})$. In particular%
\begin{equation*}
\big\|\mathcal{M(}f_{k})\big\|_{L_{p}(\mathbb{R}^{n},t_{k})}\leq c\big\|f_{k}%
\big\|_{L_{p}(\mathbb{R}^{n},t_{k})}
\end{equation*}%
holds for all sequence of functions $f_{k}\in L_{p}(\mathbb{R}^{n},t_{k})$, $%
k\in \mathbb{N}_{0}$, where $c>0$ is independent of $k$.
\end{lemma}

\begin{remark}
\label{r-estimates}We would like to mention that the result of this lemma is
true if we assume that $t_{k}\in A_{\frac{p}{\theta }}(\mathbb{R}^{n})$,\ $%
k\in \mathbb{N}_{0}$, $1<\theta \leq p<\infty $ with%
\begin{equation*}
A_{\frac{p}{\theta }}(t_{k})\leq c,\quad k\in \mathbb{N}_{0},
\end{equation*}%
where $c>0$ independent of $k$.
\end{remark}

\section{Function spaces\textbf{\ }}

In this section we\ present the Fourier analytical definition of Besov and
Triebel-Lizorkin spaces of variable smoothness and we recall their basic
properties in analogy to the classical Besov and Triebel-Lizorkin spaces. We
first need the concept of a smooth dyadic resolution of unity. Let $\varphi
_{0}$\ be a function\ in $\mathcal{S}(\mathbb{R}^{n})$\ satisfying $\varphi
_{0}(x)=1$\ for\ $\left\vert x\right\vert \leq 1$\ and\ $\varphi _{0}(x)=0$\
for\ $\left\vert x\right\vert \geq \frac{3}{2}$.\ We put $\varphi
_{k}(x)=\varphi _{0}(2^{-k}x)-\varphi _{0}(2^{1-k}x)$ for $k=1,2,3,...$.
Then $\{\varphi _{k}\}_{k\in \mathbb{N}_{0}}$\ is a resolution of unity, $%
\sum_{k=0}^{\infty }\varphi _{k}(x)=1$ for all $x\in \mathbb{R}^{n}$.\ Thus
we obtain the Littlewood-Paley decomposition 
\begin{equation*}
f=\sum_{k=0}^{\infty }\mathcal{F}^{-1}\varphi _{k}\ast f
\end{equation*}%
of all $f\in \mathcal{S}^{\prime }(\mathbb{R}^{n})$ $($convergence in $%
\mathcal{S}^{\prime }(\mathbb{R}^{n}))$.

Now, we define the spaces under consideration.

\begin{definition}
\label{B-F-def-inh}Let $0<p\leq \infty $ and $0<q\leq \infty $. Let $%
\{t_{k}\}_{k\in \mathbb{N}_{0}}$ be a $p$-admissible weight sequence. 
\newline
$\mathrm{(i)}$ The Besov space $B_{p,q}(\mathbb{R}^{n},\{t_{k}\}_{k\in 
\mathbb{N}_{0}})$\ is the collection of all $f\in \mathcal{S}^{\prime }(%
\mathbb{R}^{n})$\ such that 
\begin{equation*}
\big\|f\big\|_{B_{p,q}(\mathbb{R}^{n},\{t_{k}\}_{k\in \mathbb{N}_{0}})}:=%
\Big(\sum\limits_{k=0}^{\infty }\big\|t_{k}(\mathcal{F}^{-1}\varphi _{k}\ast
f)\big\|_{L_{p}(\mathbb{R}^{n})}^{q}\Big)^{\frac{1}{q}}<\infty
\end{equation*}%
with the usual modifications if $q=\infty $.\newline
$\mathrm{(ii)}$ Let $0<p<\infty $. The Triebel-Lizorkin space $F_{p,q}(%
\mathbb{R}^{n},\{t_{k}\}_{k\in \mathbb{N}_{0}})$\ is the collection of all $%
f\in \mathcal{S}^{\prime }(\mathbb{R}^{n})$\ such that 
\begin{equation*}
\big\|f\big\|_{F_{p,q}(\mathbb{R}^{n},\{t_{k}\}_{k\in \mathbb{N}_{0}})}:=%
\Big\|\Big(\sum\limits_{k=0}^{\infty }t_{k}^{q}|\mathcal{F}^{-1}\varphi
_{k}\ast f|^{q}\Big)^{\frac{1}{q}}\Big\|_{L_{p}(\mathbb{R}^{n})}<\infty
\end{equation*}%
with the usual modifications if $q=\infty $.
\end{definition}

Let $0<\theta \leq p<\infty $ and $0<q<\infty $. Let $\{t_{k}\}\in X_{\alpha
,\sigma ,p}$ be a $p$-admissible weight sequence with $\sigma =(\sigma
_{1}=\theta \left( \frac{p}{\theta }\right) ^{\prime },\sigma _{2}\geq p)$
and $\alpha =(\alpha _{1},\alpha _{2})\in \mathbb{R}^{2}$. The definition of
the spaces $A_{p,q}(\mathbb{R}^{n},\{t_{k}\}_{k\in \mathbb{N}_{0}})$ is
independent of the particular choice of the smooth dyadic resolution of
unity\ $\{\varphi _{k}\}_{k\in \mathbb{N}_{0}}$. They are quasi-Banach
spaces. They are Banach spaces if $1\leq p<\infty $ and $1\leq q<\infty $.\
We have the embedding%
\begin{equation*}
\mathcal{S}(\mathbb{R}^{n})\hookrightarrow A_{p,q}(\mathbb{R}%
^{n},\{t_{k}\}_{k\in \mathbb{N}_{0}})\hookrightarrow \mathcal{S}^{\prime }(%
\mathbb{R}^{n}).
\end{equation*}

Further results such as the $\varphi $-transform characterization in the
sense of Frazier and Jawerth, the case $p=\infty $, duality, complex
interpolation, the smooth atomic, molecular and wavelet decomposition and
the characterization of these function spaces in terms of the difference
relations are given in \cite{D6} and \cite{D6.1}. The above function spaces
whose elements are not distributions, but rather functions that are locally
integrable in some power are studied by Tyulenev \cite{Ty14}, \cite{Ty15}\
and \cite{Ty-N-L} for the Besov case and by the author \cite{D7} for
Triebel-Lizorkin case.

Using the system $\{\varphi _{k}\}_{k\in \mathbb{N}_{0}}$ we can define the
quasi-norms%
\begin{equation*}
\big\|f\big\|_{B_{p,q}^{s}(\mathbb{R}^{n})}:=\Big(\sum\limits_{k=0}^{\infty
}2^{ksq}\big\|\mathcal{F}^{-1}\varphi _{k}\ast f\big\|_{L_{p}(\mathbb{R}%
^{n})}^{q}\Big)^{\frac{1}{q}}
\end{equation*}%
and%
\begin{equation*}
\big\|f\big\|_{F_{p,q}^{s}(\mathbb{R}^{n})}:=\Big\|\Big(\sum\limits_{k=0}^{%
\infty }2^{ksq}|\mathcal{F}^{-1}\varphi _{k}\ast f|^{q}\Big)^{\frac{1}{q}}%
\Big\|_{L_{p}(\mathbb{R}^{n})}
\end{equation*}%
for constants $s\in \mathbb{R}$ and $0<p,q\leq \infty $ with $0<p<\infty $
in the $F$-case. The Besov space $B_{p,q}^{s}(\mathbb{R}^{n})$\ consist of
all distributions $f\in \mathcal{S}^{\prime }(\mathbb{R}^{n})$ for which $%
\big\|f\big\|_{B_{p,q}^{s}(\mathbb{R}^{n})}<\infty $. The Triebel-Lizorkin
space $F_{p,q}^{s}(\mathbb{R}^{n})$\ consist of all distributions $f\in 
\mathcal{S}^{\prime }(\mathbb{R}^{n})$ for which $\big\|f\big\|_{F_{p,q}^{s}(%
\mathbb{R}^{n})}<\infty $. It is well-known that these spaces do not depend
on the choice of the system $\{\varphi _{k}\}_{k\in \mathbb{N}_{0}}$ (up to
equivalence of quasinorms). Further details on the classical theory of these
spaces, included the homogeneous case, can be found \cite{FJ90}, \cite{T1}
and \cite{T2}.

One recognizes immediately that if $\{t_{k}\}_{k\in \mathbb{N}%
_{0}}=\{2^{sk}\}_{k\in \mathbb{N}_{0}}$, $s\in \mathbb{R}$, then we have 
\begin{equation*}
B_{p,q}(\mathbb{R}^{n},\{2^{sk}\}_{k\in \mathbb{N}_{0}})=B_{p,q}^{s}(\mathbb{%
R}^{n})
\end{equation*}%
and 
\begin{equation*}
F_{p,q}(\mathbb{R}^{n},\{2^{sk}\}_{k\in \mathbb{N}_{0}})=F_{p,q}^{s}(\mathbb{%
R}^{n}).
\end{equation*}%
Moreover, for $\{t_{k}\}_{k\in \mathbb{N}_{0}}=\{2^{sk}w\}_{k\in \mathbb{N}%
_{0}}$, $s\in \mathbb{R}$ with a weight $w$ we re-obtain the weighted Besov
and Triebel-Lizorkin spaces; we refer, in particular, to the papers \cite%
{Bui82}, \cite{BPT97} and \cite{IzSa12} for a comprehensive treatment of the
weighted spaces.

\begin{example}
A sequence $\{\gamma _{j}\}_{j\in \mathbb{N}_{0}}$ of positive real numbers
is said to be admissible if there exist two positive constants $d_{0}$ and $%
d_{1}$ such that%
\begin{equation*}
d_{0}\gamma _{j}\leq \gamma _{j+1}\leq d_{1}\gamma _{j},\quad j\in \mathbb{N}%
_{0}.
\end{equation*}%
For an admissible sequence $\{\gamma _{j}\}_{j\in \mathbb{N}_{0}}$, let%
\begin{equation*}
\underline{\gamma }_{j}=\inf_{k\geq 0}\frac{\gamma _{j+k}}{\gamma _{k}}\quad 
\text{and}\quad \overline{\gamma }_{j}=\sup_{k\geq 0}\frac{\gamma _{j+k}}{%
\gamma _{k}},\quad j\in \mathbb{N}_{0}.
\end{equation*}%
and%
\begin{equation*}
\alpha _{\gamma }=\lim_{j\longrightarrow \infty }\frac{\log \overline{\gamma 
}_{j}}{j}\quad \text{and}\quad \beta _{\gamma }=\lim_{j\longrightarrow
\infty }\frac{\log \underline{\gamma }_{j}}{j},
\end{equation*}%
be the upper and lower Boyd index of the given sequence $\{\gamma
_{j}\}_{j\in \mathbb{N}_{0}}$, respectively. Then%
\begin{equation*}
\underline{\gamma }_{j}\gamma _{k}\leq \gamma _{j+k}\leq \overline{\gamma }%
_{j}\gamma _{k},\quad j,k\in \mathbb{N}_{0}
\end{equation*}%
and for each $\varepsilon >0$,%
\begin{equation*}
c_{1}2^{(\beta _{\gamma }-\varepsilon )j}\leq \underline{\gamma }_{j}\leq 
\overline{\gamma }_{j}\leq c_{2}2^{(\alpha _{\gamma }+\varepsilon )j},\quad
j\in \mathbb{N}_{0}
\end{equation*}%
for some constants $c_{1}=c_{1}(\varepsilon )>0$ and $c_{2}=c_{2}(%
\varepsilon )>0$.\newline
Clearly the sequence $\{\gamma _{j}\}_{j\in \mathbb{N}_{0}}$ lies in $%
X_{\alpha ,\sigma ,p}$ for$\ \alpha _{1}=\beta _{\gamma }-\varepsilon
,\alpha _{2}=\alpha _{\gamma }+\varepsilon $ and $0<p,\sigma _{1},\sigma
_{2}\leq \infty $.\newline
These type of admissible sequences are used in \cite{FL06} to study Besov
and Lizorkin-Triebel spaces in terms of a generalized smoothness, see also 
\cite{HaS08}.\newline
Let us consider some examples of admissible sequences. The sequence $%
\{\gamma _{j}\}_{j\in \mathbb{N}_{0}}$,%
\begin{equation*}
\gamma _{j}=2^{sj}(1+j)^{b}(1+\log (1+j))^{c},\quad j\in \mathbb{N}_{0}
\end{equation*}%
with arbitrary fixed real numbers $s,b$ and $c$ is a an admissible sequence
with 
\begin{equation*}
\beta _{\gamma }=\alpha _{\gamma }=s.
\end{equation*}
\end{example}

\begin{example}
\label{Example1 copy(1)}Let $0<r<p<\infty $, a weight $\omega ^{p}\in A_{%
\frac{p}{r}}(\mathbb{R}^{n})$ and $\{s_{k}\}=\{2^{ks}\omega
^{p}(2^{-k})\}_{k\in \mathbb{N}_{0}}$, $s\in \mathbb{R}$. Obviously, $%
\{s_{k}\}_{k\in \mathbb{N}_{0}}$ lies in $X_{\alpha ,\sigma ,p}$ for $\alpha
_{1}=\alpha _{2}=s$, $\sigma =(r(\frac{p}{r})^{\prime },p)$.\ 
\end{example}

Let $f$ be an arbitrary function on $\mathbb{R}^{n}$ and $x,h\in \mathbb{R}%
^{n}$. Then%
\begin{equation*}
\Delta _{h}f(x):=f(x+h)-f(x),\quad \Delta _{h}^{M+1}f(x):=\Delta _{h}(\Delta
_{h}^{M}f)(x),\quad M\in \mathbb{N}.
\end{equation*}%
These are the well-known differences of functions which play an important
role in the theory of function spaces. Using mathematical induction one can
show the explicit formula%
\begin{equation*}
\Delta _{h}^{M}f(x)=\sum_{j=0}^{M}\left( -1\right) ^{j}C_{j}^{M}f(x+(M-j)h),
\end{equation*}%
where $C_{j}^{M}$ are the binomial coefficients.

Let $M\in \mathbb{N}$. For $f\in L_{1}^{\mathrm{loc}},x\in \mathbb{R}^{n}$
and a cube $Q$, we put

\begin{equation*}
\delta ^{M}(Q)f:=\frac{1}{[l(Q)]^{2n}}\int_{l(Q)I^{n}}\int_{Q}\left\vert
\Delta _{h}^{M}f(x)\right\vert dxdh,
\end{equation*}%
\begin{equation*}
\tilde{L}_{p}(\mathbb{R}^{n},t_{0}):=\Big\{f:\big\|f\big\|_{\tilde{L}_{p}(%
\mathbb{R}^{n},t_{0})}=\Big(\int_{\mathbb{R}^{n}}t_{0}^{p}(x)\big\|f\big\|%
_{L_{1}(x+I^{n})}^{p}dx\Big)^{\frac{1}{p}}<\infty \Big\},
\end{equation*}%
and%
\begin{equation*}
\delta
^{M}(x+2^{-k}I^{n})f:=2^{2kn}\int_{2^{-k}I^{n}}\int_{x+2^{-k}I^{n}}\left%
\vert \Delta _{h}^{M}f(y)\right\vert dydh,
\end{equation*}%
where $I^{n}:=(-1,1)^{n}$.

We recall the definition of the spaces $\tilde{F}_{p,q}^{M}(\mathbb{R}%
^{n},\{t_{k}\}_{k\in \mathbb{N}_{0}})$ as given in \cite{D7}.

\begin{definition}
\label{L-T-def1}Let $M\in \mathbb{N},0<p<\infty ,0<q\leq \infty $, and let $%
\{t_{k}\}_{k\in \mathbb{N}_{0}}$ be a $p$-admissible weight sequence. We set%
\begin{equation*}
\tilde{F}_{p,q}^{M}(\mathbb{R}^{n},\{t_{k}\}_{k\in \mathbb{N}_{0}}):=\Big\{%
f:f\in L_{1}^{\mathrm{loc}},\big\|f\big\|_{\tilde{F}_{p,q}^{M}(\mathbb{R}%
^{n},\{t_{k}\}_{k\in \mathbb{N}_{0}})}<\infty \Big\},
\end{equation*}%
where%
\begin{equation*}
\big\|f\big\|_{\tilde{F}_{p,q}^{M}(\mathbb{R}^{n},\{t_{k}\}_{k\in \mathbb{N}%
_{0}})}:=\big\|f\big\|_{\tilde{F}_{p,q}^{M}(\mathbb{R}^{n},\{t_{k}\}_{k\in 
\mathbb{N}_{0}})}^{\bullet }+\big\|f\big\|_{\tilde{L}_{p}(\mathbb{R}%
^{n},t_{0})},
\end{equation*}%
making the obvious modifications for $q=\infty $, with%
\begin{equation*}
\big\|f\big\|_{\tilde{F}_{p,q}^{M}(\mathbb{R}^{n},\{t_{k}\}_{k\in \mathbb{N}%
_{0}})}^{\bullet }:=\Big\|\Big(\sum\limits_{k=1}^{\infty }t_{k}^{q}\left(
\delta ^{M}(\cdot +2^{-k}I^{n})f\right) ^{q}\Big)^{\frac{1}{q}}\Big\|_{L_{p}(%
\mathbb{R}^{n})}.
\end{equation*}
\end{definition}

Now we present the definition of Besov spaces of variable smoothness $\tilde{%
B}_{p,q}^{M}(\mathbb{R}^{n},\{t_{k}\}_{k\in \mathbb{N}_{0}})$ as introduced
recently in \cite{Ty15}.

\begin{definition}
Let $M\in \mathbb{N},0<p,q\leq \infty $, and let $\{t_{k}\}_{k\in \mathbb{N}%
_{0}}$ be a $p$-admissible weight sequence. We set%
\begin{equation*}
\tilde{B}_{p,q}^{M}(\mathbb{R}^{n},\{t_{k}\}_{k\in \mathbb{N}_{0}}):=\Big\{%
f:f\in L_{1}^{\mathrm{loc}},\big\|f\big\|_{\tilde{B}_{p,q}^{M}(\mathbb{R}%
^{n},\{t_{k}\}_{k\in \mathbb{N}_{0}})}<\infty \Big\},
\end{equation*}%
where%
\begin{equation*}
\big\|f\big\|_{\tilde{B}_{p,q}^{M}(\mathbb{R}^{n},\{t_{k}\}_{k\in \mathbb{N}%
_{0}}}:=\Big(\sum\limits_{k=1}^{\infty }\big\|t_{k}\delta ^{M}(\cdot
+2^{-k}I^{n})f\big\|_{L_{p}(\mathbb{R}^{n})}^{q}\Big)^{\frac{1}{q}}+\big\|f%
\big\|_{\tilde{L}_{p}(\mathbb{R}^{n},t_{0})},
\end{equation*}%
making the obvious modifications for $p=\infty $ and/or $q=\infty $.
\end{definition}

Let $M\in \mathbb{N},0<p<\infty ,0<q\leq \infty $, and let $\{t_{k}\}_{k\in 
\mathbb{N}_{0}}$ be a $p$-admissible weight sequence. We set%
\begin{equation*}
\big\|f\big\|_{\tilde{F}_{p,q}^{M}(\mathbb{R}^{n},\{t_{k}\}_{k\in \mathbb{N}%
_{0}})}^{\ast }:=\big\|f\big\|_{\tilde{F}_{p,q}^{M}(\mathbb{R}%
^{n},\{t_{k}\}_{k\in \mathbb{N}_{0}})}^{\ast ,1}+\Big(\sum_{m\in \mathbb{Z}%
^{n}}t_{0,m}^{p}\big\|f\big\|_{L_{1}(Q_{0,m})}^{p}\Big)^{\frac{1}{p}},
\end{equation*}%
where 
\begin{equation*}
\big\|f\big\|_{\tilde{F}_{p,q}^{M}(\mathbb{R}^{n},\{t_{k}\}_{k\in \mathbb{N}%
_{0}})}^{\ast ,1}:=\Big\|\Big(\sum\limits_{k=1}^{\infty }\sum_{m\in \mathbb{Z%
}^{n}}2^{kn\frac{q}{p}}t_{k,m}^{q}(\delta ^{M}(Q_{k,\tilde{m}})f)^{q}\chi
_{k,m}\Big)^{\frac{1}{q}}\Big\|_{L_{p}(\mathbb{R}^{n})},
\end{equation*}%
Also we set, $0<p\leq \infty ,$%
\begin{equation*}
\big\|f\big\|_{\tilde{B}_{p,q}^{M}(\mathbb{R}^{n},\{t_{k}\}_{k\in \mathbb{N}%
_{0}})}^{\ast }:=\big\|f\big\|_{\tilde{B}_{p,q}^{M}(\mathbb{R}%
^{n},\{t_{k}\}_{k\in \mathbb{N}_{0}})}^{\ast ,1}+\Big(\sum_{m\in \mathbb{Z}%
^{n}}t_{0,m}^{p}\big\|f\big\|_{L_{1}(Q_{0,m})}^{p}\Big)^{\frac{1}{p}},
\end{equation*}%
where 
\begin{equation*}
\big\|f\big\|_{\tilde{B}_{p,q}^{M}(\mathbb{R}^{n},\{t_{k}\}_{k\in \mathbb{N}%
_{0}})}^{\ast ,1}:=\Big(\sum\limits_{k=1}^{\infty }\Big(\sum_{m\in \mathbb{Z}%
^{n}}t_{k,m}^{p}(\delta ^{M}(Q_{k,m})f)^{p}\Big)^{\frac{q}{p}}\Big)^{\frac{1%
}{q}},
\end{equation*}%
with%
\begin{equation*}
Q_{k,\tilde{m}}:=\prod_{i=1}^{n}\big(\frac{m_{i}-2}{2^{k}},\frac{m_{i}+3}{%
2^{k}}\big),\quad \chi _{k,m}:=\chi _{Q_{k,m}},\quad m\in \mathbb{Z}^{n}
\end{equation*}%
and 
\begin{equation*}
\delta ^{M}(Q_{k,\tilde{m}})f:=\frac{1}{[l(Q_{k,\tilde{m}}]^{2n}}%
\int_{2^{-k}I^{n}}\int_{Q_{k,\tilde{m}}}\left\vert \Delta
_{h}^{M}f(z)\right\vert dzdh.
\end{equation*}

For simplicity, in what follows, we use the notation $\tilde{A}_{p,q}^{M}(%
\mathbb{R}^{n},\{t_{k}\}_{k\in \mathbb{N}_{0}})$ to denote either $\tilde{B}%
_{p,q}^{M}(\mathbb{R}^{n},\{t_{k}\}_{k\in \mathbb{N}_{0}})$ or $\tilde{F}%
_{p,q}^{M}(\mathbb{R}^{n},\{t_{k}\}_{k\in \mathbb{N}_{0}})$. The following
theorems are useful for us.

\begin{theorem}
\label{Equi-norm}Let $\alpha _{1},\alpha _{2}\in \mathbb{R}$ and $\alpha
=(\alpha _{1},\alpha _{2})$. Let $M\in \mathbb{N},0<\theta \leq p<\infty \ $%
and $0<q<\infty $. Let $\{t_{k}\}_{k\in \mathbb{N}_{0}}\in X_{\alpha ,\sigma
,p}$ be a $p$-admissible sequence with $\sigma =(\sigma _{1},\sigma _{2})$, $%
\sigma _{1}=\theta \left( \frac{p}{\theta }\right) ^{\prime }$ and $\sigma
_{2}\geq p$. Then%
\begin{equation*}
\big\|\cdot \big\|_{\tilde{A}_{p,q}^{M}(\mathbb{R}^{n},\{t_{k}\}_{k\in 
\mathbb{N}_{0}})}^{\ast }
\end{equation*}%
is an equivalent quasi-norm in $\tilde{A}_{p,q}^{M}(\mathbb{R}%
^{n},\{t_{k}\}_{k\in \mathbb{N}_{0}})$.
\end{theorem}

This theorem for Besov case is given in \cite{Ty15}, while the
Triebel-Lizorkin case can be proved as in \cite{D7}.

In the following theorem we present the characterizations of $A_{p,q}(%
\mathbb{R}^{n},\{t_{k}\}_{k\in \mathbb{N}_{0}})$ in terms of the difference
relations, see \cite{D6} for Besov case, while for Triebel-Lizorkin spaces
can be obtained by the similar arguments.

\begin{theorem}
\label{means-diff-cha1}\textit{Let }$1\leq \theta <p<\infty ,1\leq q<\infty
,\alpha _{1},\alpha _{2}\in \mathbb{R},\alpha =(\alpha _{1},\alpha _{2})$
and $M\in \mathbb{N}$. Let $\{t_{k}\}_{k\in \mathbb{N}_{0}}\in X_{\alpha
,\sigma ,p}$ be a $p$-admissible weight sequence with $\sigma =(\sigma
_{1}=\theta \left( \frac{p}{\theta }\right) ^{\prime },\sigma _{2}\geq p)$.
Assume that%
\begin{equation*}
0<\alpha _{1}\leq \alpha _{2}<M.
\end{equation*}%
Then 
\begin{equation*}
B_{p,q}(\mathbb{R}^{n},\{t_{k}\}_{k\in \mathbb{N}_{0}})=\tilde{B}_{p,q}^{M}(%
\mathbb{R}^{n},\{t_{k}\}_{k\in \mathbb{N}_{0}})
\end{equation*}%
and 
\begin{equation*}
F_{p,q}(\mathbb{R}^{n},\{t_{k}\}_{k\in \mathbb{N}_{0}})=\tilde{F}_{p,q}^{M}(%
\mathbb{R}^{n},\{t_{k}\}_{k\in \mathbb{N}_{0}}),
\end{equation*}%
in the sense of equivalent norm.
\end{theorem}

\section{Proof of the main result}

First, we prove that our estimate $\mathrm{\eqref{key-est}}$ is better than $%
\mathrm{\eqref{key-est1}.}$ Obviously 
\begin{equation*}
H\leq \sup_{k\geq i,x\in \mathbb{R}^{n}}\frac{t_{k-i}(\lambda ^{-1}x)}{%
t_{k-i}(x)}.
\end{equation*}%
Let $\{t_{k}\}_{k\in \mathbb{N}_{0}}=\{2^{ks}\omega \}_{k\in \mathbb{N}_{0}}$%
, with $\omega (x)=|x-1|^{\delta }$, $s>0$ and $\delta ,x\in \mathbb{R}$.
Then $t_{k}^{p}\in A_{\frac{p}{\theta }}(\mathbb{R})$, $1<\frac{p}{\theta }%
<\infty $, if and only if $-\frac{1}{p}<\delta <\frac{1}{\theta }-\frac{1}{p}
$. We have%
\begin{equation*}
\sup_{x\in \mathbb{R}^{n}}\frac{\omega (\lambda ^{-1}x)}{\omega (x)}\geq
\lambda ^{-\delta }\sup_{|x|\geq \lambda }\frac{|x-\lambda |^{\delta }}{%
|x-1|^{\delta }}\geq \lambda ^{-\delta }\sup_{|x|\geq \lambda }|x-\lambda
|^{\delta }=\infty
\end{equation*}%
for any $-\frac{1}{p}<\delta <0$ and any $\lambda >1$. Let $p_{0}=\frac{p}{%
\theta }$. From Lemma \ref{Ap-Property}/(iv) we conclude%
\begin{equation*}
\frac{M_{Q_{k-i,m},p}(\omega (\lambda ^{-1}\cdot ))}{M_{Q_{k-i,m},p}(\omega )%
}\lesssim \frac{\Big(M_{Q_{k-i,m},\theta p_{0}^{\prime }}\big(\left( \omega
(\lambda ^{-1}\cdot )\right) ^{-1}\big)\Big)^{-1}}{M_{Q_{k-i,m},p}(\omega )}.
\end{equation*}%
We have%
\begin{eqnarray*}
M_{Q_{k-i,m},\theta p_{0}^{\prime }}\big(\left( \omega (\lambda ^{-1}\cdot
)\right) ^{-1}\big) &\geq &M_{Q_{k-i,m},\theta p_{0}^{\prime }}\big(\left(
\omega (\lambda ^{-1}\cdot )\right) ^{-1}\chi _{|\cdot -1|\geq 2(\lambda -1)}%
\big) \\
&\geq &cM_{Q_{k-i,m},\theta p_{0}^{\prime }}(\omega ^{-1}\chi _{|\cdot
-1|\geq 2(\lambda -1)}),
\end{eqnarray*}%
which yields that%
\begin{equation*}
\frac{M_{Q_{k-i,m},p}(\omega (\lambda ^{-1}\cdot ))}{M_{Q_{k-i,m},p}(\omega )%
}\lesssim \frac{\left( M_{Q_{k-i,m},\theta p_{0}^{\prime }}(\omega ^{-1}\chi
_{|\cdot -1|\geq 2(\lambda -1)})\right) ^{-1}}{M_{Q_{k-i,m},p}(\omega \chi
_{|\cdot -1|\geq 2(\lambda -1)})}\leq c,\quad m\in \mathbb{Z}\text{, }k\geq i
\end{equation*}%
because of $\omega ^{p}\chi _{|\cdot -1|\geq 2(\lambda -1)}\in A_{\frac{p}{%
\theta }}(\mathbb{R})$, where $c$ independent of $k,i$ and $m$. Therefore $H$
is finite.

Now we prove Theorem \ref{Dilation}. By similarity, we only consider the
space $F_{p,q}(\mathbb{R}^{n},\{t_{k}\}_{k\in \mathbb{N}_{0}})$. Let $M\in 
\mathbb{N}$ be such that $0<\alpha _{1}\leq \alpha _{2}<M$. Of course, $%
f(\lambda \cdot )$ must be interpreted in the sense of distributions. On the
other hand by the embedding 
\begin{equation*}
F_{p,q}(\mathbb{R}^{n},\{t_{k}\}_{k\in \mathbb{N}_{0}})\hookrightarrow
L_{1}^{\mathrm{loc}},
\end{equation*}%
if $\alpha _{2}\geq \alpha _{1}>0$, and $1\leq p,q<\infty $, see \cite{D6},
it follows that $f(x)$ is a regular distribution and $f(\lambda x)$ makes
also sense as a locally integrable function. Since $0<\alpha _{1}\leq \alpha
_{2}<M$, by Theorem \ref{means-diff-cha1} 
\begin{equation*}
\big\|f\big\|_{\tilde{F}_{p,q}^{M}(\mathbb{R}^{n},\{t_{k}\}_{k\in \mathbb{N}%
_{0}})}
\end{equation*}%
is an equivalent norm in $F_{p,q}(\mathbb{R}^{n},\{t_{k}\}_{k\in \mathbb{N}%
_{0}})$. We will prove 
\begin{equation*}
\big\|f(\lambda \cdot )\big\|_{\tilde{F}_{p,q}^{M}(\mathbb{R}%
^{n},\{t_{k}\}_{k\in \mathbb{N}_{0}})}\leq c\lambda ^{\alpha _{2}-\frac{n}{p}%
}H\big\|f|\big\|_{\tilde{F}_{p,q}^{M}(\mathbb{R}^{n},\{t_{k}\}_{k\in \mathbb{%
N}_{0}})}
\end{equation*}%
for all $0<p,q<\infty $ and 
\begin{equation*}
\sigma _{p}=\max \Big(0,\frac{n}{p}-n\Big)<\alpha _{1}\leq \alpha _{2}<M.
\end{equation*}

\textbf{Step 1.} After a simple change of variable%
\begin{equation*}
\big\|f(\lambda \cdot )\big\|_{\tilde{L}_{p}(\mathbb{R}^{n},t_{0})}\leq
\lambda ^{-\frac{n}{p}}\Big(\lambda ^{-np}\int_{\mathbb{R}%
^{n}}t_{0}^{p}(\lambda ^{-1}x)\big\|f\big\|_{L_{1}(x+2^{i}I^{n})}^{p}dx\Big)%
^{\frac{1}{p}}.
\end{equation*}%
The above expression in the brackets can be rewritten as follows%
\begin{equation*}
\lambda ^{-np}\sum_{m\in \mathbb{Z}^{n}}\int_{Q_{-i,m}}t_{0}^{p}(\lambda
^{-1}x)\big\|f\big\|_{L_{1}(x+2^{i}I^{n})}^{p}dx,
\end{equation*}%
which can be estimated\ by%
\begin{equation}
\lambda ^{-np}\sum_{m\in \mathbb{Z}^{n}}\big\|f\big\|_{L_{1}(\tilde{Q}%
_{-i,m})}^{p}\big\|t_{0}(\lambda ^{-1}\cdot )\big\|_{L_{p}(Q_{-i,m})}^{p},
\label{mainterm}
\end{equation}%
where $\tilde{Q}_{-i,m}=\cup _{j=1}^{3^{n}}Q_{-i,z_{j}(m)}$. Therefore $%
\mathrm{\eqref{mainterm}}$ can be rewritten in the following form%
\begin{equation*}
\lambda ^{-np}\sum_{j=1}^{3^{n}}\sum_{m\in \mathbb{Z}^{n}}\big\|%
t_{0}(\lambda ^{-1}\cdot )\big\|_{L_{p}(Q_{-i,m})}^{p}\Big(\sum_{h\in 
\mathbb{Z}^{n},Q_{0,h}\subset Q_{-i,z_{j}(m)}}\big\|f\big\|_{L_{1}(Q_{0,h})}%
\Big)^{p},
\end{equation*}%
which is bounded by%
\begin{equation}
c\lambda ^{np(\frac{1}{\theta }-1)}H^{p}\sum_{j=1}^{3^{n}}\sum_{m\in \mathbb{%
Z}^{n}}\Big(\sum_{h\in \mathbb{Z}^{n},Q_{0,h}\subset Q_{-i,z_{j}(m)}}t_{0,h}%
\big\|f\big\|_{L_{1}(Q_{0,h})}\Big)^{p},  \label{key-term}
\end{equation}%
where we have used Lemma \ref{Ap-Property}/(iii) and the constant $c$ is
independent of $\lambda $. We distinguish two cases.

\textbf{Case 1.} $1<p<\infty $. Observe that $\left\vert
h-2^{i}z_{j}(m)\right\vert \leq c2^{i}$ for any $h\in \mathbb{Z}^{n}$, such
that $Q_{0,h}\subset Q_{k-i,z_{j}(m)}$, $j=1,...,3^{n}$. Therefore the
number of terms in the sum $\sum_{h\in \mathbb{Z}^{n}:Q_{0,h}\subset
Q_{k-i,z_{j}(m)}}\cdot \cdot \cdot $ does not exceed $c2^{in}$. By H\"{o}%
lder's inequality $\mathrm{\eqref{key-term}}$ is bounded by 
\begin{equation*}
c\lambda ^{n(\frac{p}{\theta }-1)}H^{p}\sum_{j=1}^{3^{n}}\sum_{m\in \mathbb{Z%
}^{n}}\sum_{h\in \mathbb{Z}^{n},Q_{0,h}\subset Q_{-i,z_{j}(m)}}t_{0,h}^{p}%
\big\|f\big\|_{L_{1}(Q_{0,h})}^{p},
\end{equation*}%
which is bounded by%
\begin{eqnarray*}
&&c\lambda ^{n(\frac{p}{\theta }-1)}H^{p}\sum_{j=1}^{3^{n}}\sum_{m\in 
\mathbb{Z}^{n}}\int_{Q_{-i,z_{j}(m)}}\sum_{h\in \mathbb{Z}^{n}}t_{0,h}^{p}%
\big\|f\big\|_{L_{1}(Q_{0,h})}^{p}\chi _{Q_{0,h}}(x)dx \\
&\lesssim &\lambda ^{n(\frac{p}{\theta }-1)}H^{p}\int_{\mathbb{R}%
^{n}}\sum_{h\in \mathbb{Z}^{n}}t_{0,h}^{p}\big\|f\big\|_{L_{1}(Q_{0,h})}^{p}%
\chi _{Q_{0,h}}(x)dx \\
&\lesssim &\lambda ^{n(\frac{p}{\theta }-1)}H^{p}\sum_{h\in \mathbb{Z}%
^{n}}t_{0,h}^{p}\big\|f\big\|_{L_{1}(Q_{0,h})}^{p}.
\end{eqnarray*}

\textbf{Case 2.} $0<p\leq 1$. We have $\mathrm{\eqref{key-term}}$ is bounded
by%
\begin{eqnarray*}
&&c\lambda ^{np(\frac{1}{\theta }-1)}H^{p}\sum_{m\in \mathbb{Z}%
^{n}}\sum_{h\in \mathbb{Z}^{n},Q_{0,h}\subset \tilde{Q}_{-i,m}}%
\int_{Q_{0,h}}t_{0,h}^{p}(x)\big\|f\big\|_{L_{1}(Q_{0,h})}^{p}dx \\
&\lesssim &\lambda ^{np(\frac{1}{\theta }-1)}H^{p}\sum_{h\in \mathbb{Z}%
^{n}}t_{0,h}^{p}\big\|f\big\|_{L_{1}(Q_{0,h})}^{p}.
\end{eqnarray*}%
Taking $0<\theta <p<\infty $ be such that $\theta >\frac{n}{\alpha _{2}+%
\frac{n}{\max (1,p)}}$\ which is possible because of $\alpha _{2}>\sigma
_{p} $. Therefore 
\begin{equation*}
\big\|f(\lambda \cdot )\big\|_{\tilde{L}_{p}(\mathbb{R}^{n},t_{0})}\leq
c\lambda ^{\alpha _{2}-\frac{n}{p}}H\big\|f|\big\|_{F_{p,q}(\mathbb{R}%
^{n},\{t_{k}\}_{k\in \mathbb{N}_{0}})}
\end{equation*}%
for some positive constant $c$ independent of $\lambda $.

\textbf{Step 2.} We will estimate 
\begin{equation}
\Big\|\Big(\sum\limits_{k=1}^{\infty }t_{k}^{q}(\delta
^{M}(x+2^{-k}I^{n})f(\lambda \cdot ))^{q}\Big)^{\frac{1}{q}}\Big\|_{L_{p}(%
\mathbb{R}^{n})}.  \label{sec-est2}
\end{equation}%
It is easily seen that%
\begin{equation*}
\delta ^{M}(x+2^{-k}I^{n})f(\lambda \cdot )\leq 2^{2n}\delta ^{M}(\lambda
x+2^{i-k}I^{n})f
\end{equation*}%
for any $x\in \mathbb{R}^{n}$. Therefore \eqref{sec-est2} can be estimated
by 
\begin{equation*}
\lambda ^{-\frac{n}{p}}\Big\|\Big(\sum\limits_{k=1}^{\infty }\big(%
t_{k}(\lambda ^{-1}\cdot )\delta ^{M}(\cdot +2^{i-k}I^{n})f\big)^{q}\Big)^{%
\frac{1}{q}}\Big\|_{L_{p}(\mathbb{R}^{n})}.
\end{equation*}%
Obviously%
\begin{eqnarray*}
&&\sum\limits_{k=1}^{\infty }\big(t_{k}(\lambda ^{-1}x)\delta
^{M}(x+2^{i-k}I^{n})f\big)^{q} \\
&=&\sum\limits_{k=1}^{\infty }\sum_{m\in \mathbb{Z}^{n}}\big(t_{k}(\lambda
^{-1}x)\delta ^{M}(x+2^{i-k}I^{n})f\big)^{q}\chi _{Q_{k-i,m}}(x).
\end{eqnarray*}%
Let $x\in Q_{k-i,m}$ with $k\in \mathbb{N}$ and $m\in \mathbb{Z}^{n}$. We
find\ that%
\begin{eqnarray*}
\delta ^{M}(x+2^{i-k}I^{n})(f) &\leq &\delta ^{M}(Q_{k-i,m}+2^{i-k}I^{n})f \\
&=&\Big(\frac{1}{|Q_{k-i,m}|}\int_{Q_{k-i,m}}[\delta
^{M}(Q_{k-i,m}+2^{i-k}I^{n})f]^{\delta }dy\Big)^{\frac{1}{\delta }},
\end{eqnarray*}%
where $0<\delta <\min (1,\theta ,q)$. Observing that 
\begin{equation*}
Q_{k-i,m}+2^{i-k}I^{n}=y+Q_{k-i,m}-y+2^{i-k}I^{n},\quad y\in Q_{k-i,m}.
\end{equation*}%
We have $Q_{k-i,m}-y\subset 2^{i-k}I^{n}$\ for all $y\in Q_{k-i,m}$ and this
implies that%
\begin{equation*}
Q_{k-i,m}+2^{i-k}I^{n}\subset y+2^{i-k+1}I^{n},\quad y\in Q_{k-i,m}.
\end{equation*}%
Therefore, for any $x\in Q_{k-i,m}$,%
\begin{equation*}
\delta ^{M}(x+2^{i-k}I^{n})(f)\leq c\big(M_{Q_{k-i,m}}([\bar{\delta}%
^{M}(\cdot +2^{i-k+1}I^{n})f]^{\delta })\big)^{\frac{1}{\delta }},
\end{equation*}%
where $c>0$ is independent of $k,m$ and $x$, with 
\begin{equation*}
\bar{\delta}^{M}(y+2^{i-k+1}I^{n})f=2^{2(k-i)n}\int_{2^{i-k}I^{n}}%
\int_{y+2^{i-k+1}I^{n}}\left\vert \Delta _{h}^{M}f(\nu )\right\vert d\nu dh.
\end{equation*}%
Hence%
\begin{equation*}
\sum\limits_{k=1}^{\infty }\big(t_{k}(\lambda ^{-1}x)\delta
^{M}(x+2^{i-k}I^{n})f\big)^{q}
\end{equation*}%
can be estimated by%
\begin{eqnarray*}
&&c\sum\limits_{k=1}^{\infty }\sum_{m\in \mathbb{Z}^{n}}t_{k}^{q}(\lambda
^{-1}x)\big(M_{Q_{k-i,m}}([\bar{\delta}^{M}(\cdot +2^{i-k+1}I^{n})f]^{\delta
})\big)^{\frac{q}{\delta }}\chi _{k-i,m}(x) \\
&=&c\sum\limits_{k=1}^{\infty }\Big(\sum_{m\in \mathbb{Z}^{n}}t_{k}^{\delta
}(\lambda ^{-1}x)M_{Q_{k-i,m}}([\bar{\delta}^{M}(\cdot
+2^{i-k+1}I^{n})f]^{\delta })\chi _{k-i,m}(x)\Big)^{\frac{q}{\delta }}
\end{eqnarray*}%
for all $x\in \mathbb{R}^{n}$, where $t_{k}^{q}(\lambda ^{-1}\cdot )=\left(
t_{k}(\lambda ^{-1}\cdot )\right) ^{q}$. Using this estimate, the quantity 
\begin{equation*}
\Big\|\Big(\sum_{k=1}^{\infty }t_{k}^{q}(\delta ^{M}(\cdot
+2^{-k}I^{n})f(\lambda \cdot ))^{q}\Big)^{1/q}\Big\|_{L_{p}(\mathbb{R}^{n})}
\end{equation*}%
can be estimated from above by%
\begin{equation*}
c\Big\|\Big(\sum\limits_{k=1}^{\infty }\Big(\sum_{m\in \mathbb{Z}%
^{n}}t_{k}^{\delta }(\lambda ^{-1}\cdot )M_{Q_{k-i,m}}([\bar{\delta}%
^{M}(\cdot +2^{i-k+1}I^{n})f]^{\delta })\chi _{k-i,m}\Big)^{\frac{q}{\delta }%
}\Big)^{\frac{\delta }{q}}\Big\|_{L_{\frac{p}{\delta }}(\mathbb{R}^{n})}^{%
\frac{1}{\delta }}.
\end{equation*}%
By duality the last term with power $\delta $ is bounded by 
\begin{equation*}
c\sup \sum\limits_{k=1}^{\infty }\sum_{m\in \mathbb{Z}^{n}}%
\int_{Q_{k-i,m}}t_{k}^{\delta }(\lambda ^{-1}x)M_{Q_{k-i,m}}([\bar{\delta}%
^{M}(\cdot +2^{i-k+1}I^{n})f]^{\delta })|g_{k}(x)|dx,
\end{equation*}%
where the supremum is taking over all sequence of functions $\{g_{k}\}_{k\in 
\mathbb{N}}\in L_{(\frac{p}{\delta })^{\prime }}(\ell _{(\frac{q}{\delta }%
)^{\prime }})$ with 
\begin{equation*}
\big\|\{g_{k}\}_{k\in \mathbb{N}}\big\|_{L_{(\frac{p}{\delta })^{\prime
}}(\ell _{(\frac{q}{\delta })^{\prime }})}\leq 1.
\end{equation*}%
By H\"{o}lder's inequality,%
\begin{equation*}
1=M_{Q_{k-i,m},h}\big(t_{k}^{-\delta }(\lambda ^{-1}\cdot )t_{k}^{\delta
}(\lambda ^{-1}\cdot )\big)\leq M_{Q_{k-i,m},\tau }\big(t_{k}^{-\delta
}(\lambda ^{-1}\cdot )\big)\big(M_{Q_{k-i,m},p}(t_{k}^{\delta }(\lambda
^{-1}\cdot )\big)
\end{equation*}%
for any $h,\tau >0$ with\ $\frac{1}{h}=\frac{1}{p}+\frac{1}{\tau }$ and $%
m\in \mathbb{Z}^{n},k\in \mathbb{N}$. By Jensen's inequality, the second
term is bounded by%
\begin{equation*}
\Big(M_{Q_{k-i,m},p}\big(t_{k}(\lambda ^{-1}\cdot )\big)\Big)^{\delta }.
\end{equation*}%
Using the fact that%
\begin{equation*}
M_{Q_{k-i,m},\tau }\big(t_{k}^{-\delta }(\lambda ^{-1}\cdot )\big)%
M_{Q_{k-i,m}}(t_{k}^{\delta }(\lambda ^{-1}\cdot )g_{k})\lesssim \big(%
M_{Q_{k-i,m}}((t_{k}^{-\delta }(\lambda ^{-1}\cdot )\mathcal{M}%
(t_{k}^{\delta }(\lambda ^{-1}\cdot )g_{k}))^{\tau }\big)^{\frac{1}{\tau }},
\end{equation*}%
we find that%
\begin{equation*}
M_{Q_{k-i,m}}([\bar{\delta}^{M}(\cdot +2^{i-k+1}I^{n})f]^{\delta
})\int_{Q_{k-i,m}}t_{k}^{\delta }(\lambda ^{-1}x)|g_{k}(x)|dx
\end{equation*}%
can be estimated by%
\begin{equation*}
c|Q_{k-i,m}|^{-\frac{\delta }{p}}\int_{Q_{k-i,m}}M_{Q_{k-i-1,m}}(\omega
_{k,i,m})\big(M_{Q_{k-i-1,m}}((t_{k}^{-\delta }(\lambda ^{-1}\cdot )\mathcal{%
M}(t_{k}^{\delta }(\lambda ^{-1}\cdot )g_{k}))^{\tau }\big)^{\frac{1}{\tau }%
}dx,
\end{equation*}%
where%
\begin{equation*}
\omega _{k,i,m}=(\tilde{t}_{k-i,m}\bar{\delta}^{M}(\cdot
+2^{i-k+1}I^{n})f)^{\delta }\quad \text{and}\quad \tilde{t}_{k-i,m}=\big\|%
t_{k}(\lambda ^{-1}\cdot )\big\|_{L_{p}(Q_{k-i,m})},
\end{equation*}%
which is bounded by%
\begin{equation*}
c|Q_{k-i,m}|^{-\frac{\delta }{p}}\int_{Q_{k-i,m}}\mathcal{M}\big(\omega
_{k,i,m}\chi _{k-i,m}\big)(x)\mathcal{M}_{\tau }\big(t_{k}^{-\delta
}(\lambda ^{-1}\cdot )\mathcal{M}(t_{k}^{\delta }(\lambda ^{-1}\cdot )g_{k})%
\big)(x)dx.
\end{equation*}%
Therefore,%
\begin{eqnarray*}
&&\sum\limits_{k=1}^{\infty }\int \sum_{m\in \mathbb{Z}^{n}}M_{Q_{k-i,m}}([%
\bar{\delta}^{M}(\cdot +2^{i-k+1}I^{n})f]^{\delta })t_{k}^{\delta }(\lambda
^{-1}x)|g_{k}(x)|\chi _{k-i,m}(x)dx \\
&\lesssim &\int \sum\limits_{k=1}^{\infty }\mathcal{M}\big(\sum_{z\in 
\mathbb{Z}^{n}}2^{(k-i)\frac{n}{p}\delta }\omega _{k,i,z}\chi _{k-i,z}\big)%
(x)\mathcal{M}_{\tau }\big(t_{k}^{-\delta }(\lambda ^{-1}\cdot )\mathcal{M}%
(t_{k}^{\delta }(\lambda ^{-1}\cdot )g_{k})\big)(x)dx.
\end{eqnarray*}%
By H\"{o}lder's inequality the term inside the integral is bounded by%
\begin{equation*}
\Big(\sum\limits_{k=1}^{\infty }\Big(\mathcal{M}\big(\sum_{z\in \mathbb{Z}%
^{n}}2^{(k-i)\frac{n}{p}\delta }\omega _{k,i,z}\chi _{k-i,z}\big)\Big)^{%
\frac{q}{\delta }}\Big)^{\frac{\delta }{q}}\Big(\sum\limits_{k=1}^{\infty }%
\big(\mathcal{M}_{\tau }(t_{k}^{-\delta }(\lambda ^{-1}\cdot )\mathcal{M}%
(t_{k}^{\delta }(\lambda ^{-1}\cdot )g_{k})\big)^{(\frac{q}{\delta }%
)^{\prime }}\Big)^{\frac{1}{(\frac{q}{\delta })^{\prime }}}.
\end{equation*}%
Since $t_{k}^{p}\in A_{\frac{p}{\theta }}(\mathbb{R}^{n})\subset A_{\frac{p}{%
\delta }}(\mathbb{R}^{n})$, using Lemma \ref{Ap-Property}/(ii), (iv) to
obtain 
\begin{equation*}
t_{k}^{-\delta (\frac{p}{\delta })^{\prime }}\in A_{(\frac{p}{\delta }%
)^{\prime }}(\mathbb{R}^{n}),\quad k\in \mathbb{N}_{0}
\end{equation*}%
and there exist a $1<\varrho <(\frac{p}{\delta })^{\prime }<\infty $ such
that 
\begin{equation*}
t_{k}^{-\delta (\frac{p}{\delta })^{\prime }}\in A_{\frac{(\frac{p}{\delta }%
)^{\prime }}{\varrho }}(\mathbb{R}^{n}),\quad k\in \mathbb{N}_{0}.
\end{equation*}%
Taking any $0<\tau <\min (1,(\frac{q}{\delta })^{\prime },(\frac{p}{\delta }%
)^{\prime })$, using the vector-valued maximal inequality of Fefferman and
Stein $\mathrm{\eqref{Fe-St71}}$, and Lemma \ref{key-estimate1} we find that
the second term of the last expression in $L_{(\frac{p}{\delta })^{\prime }}(%
\mathbb{R}^{n})$-norm is bounded, while the first term in $L_{\frac{p}{%
\delta }}(\mathbb{R}^{n})$-norm is bounded by 
\begin{eqnarray}
&&c\Big\|\Big(\sum\limits_{k=1}^{\infty }\Big(\mathcal{M}\big(2^{(k-i)\frac{n%
}{p}\delta }\sum_{z\in \mathbb{Z}^{n}}\omega _{k,i,z}\chi _{k-i,z}\big)\Big)%
^{\frac{q}{\delta }}\Big)^{\frac{\delta }{q}}\Big\|_{L_{\frac{p}{\delta }}(%
\mathbb{R}^{n})}  \notag \\
&\lesssim &\Big\|\Big(\sum\limits_{k=1}^{\infty }\sum_{z\in \mathbb{Z}%
^{n}}2^{(k-i)\frac{n}{p}q}\big(\tilde{t}_{k-i,z}\bar{\delta}^{M}(\cdot
+2^{i-k+2}I^{n})f\big)^{q}\chi _{k-i,z}\Big)^{\frac{1}{q}}\Big\|_{L_{p}(%
\mathbb{R}^{n})}^{\delta },  \label{est3}
\end{eqnarray}%
where we used the vector-valued maximal inequality of Fefferman and Stein $%
\mathrm{\eqref{Fe-St71}}$ since $0<\delta <\min (1,\theta ,q)$. Now by $%
\mathrm{\eqref{Asum2}}$ we deduce that 
\begin{equation*}
\tilde{t}_{k-i,m}\lesssim 2^{\alpha _{2}i}\big\|t_{k-i}(\lambda ^{-1}\cdot )%
\big\|_{L_{p}(Q_{k-i,m})},\quad m\in \mathbb{Z}^{n},k\geq i.
\end{equation*}%
Therefore $\mathrm{\eqref{est3}}$, with $\sum\limits_{k=i+1}^{\infty }$ in
place of $\sum\limits_{k=1}^{\infty }$ is bounded by%
\begin{equation}
c\lambda ^{\alpha _{2}\delta }H^{\delta }\Big\|\Big(\sum\limits_{j=1}^{%
\infty }\sum_{z\in \mathbb{Z}^{n}}2^{j\frac{n}{p}}\big(t_{j,z}\bar{\delta}%
^{M}(\cdot +2^{1-j}I^{n})f\big)^{q}\chi _{j,z}\Big)^{\frac{1}{q}}\Big\|%
_{L_{p}(\mathbb{R}^{n})}^{\delta }.  \label{est-norm}
\end{equation}%
Using the fact that%
\begin{equation*}
\bar{\delta}^{M}(\cdot +2^{1-j}I^{n})f\leq \delta ^{M}(Q_{j,\tilde{z}})f
\end{equation*}%
and Theorem \ref{Equi-norm} we easily estimate $\mathrm{\eqref{est-norm}}$ by%
\begin{equation*}
c\lambda ^{\alpha _{2}\delta }H^{\delta }\big\Vert f\big\Vert_{F_{p,q}(%
\mathbb{R}^{n},\{t_{k}\}_{k\in \mathbb{N}_{0}})}^{\delta },
\end{equation*}%
where the constant $c$ independent of $\lambda $.

Now for any $x\in Q_{k-i,z}$, with $k\in \{1,...,i\}$ and $z\in \mathbb{Z}%
^{n}$ we have%
\begin{equation*}
\bar{\delta}^{M}(x+2^{i-k+1}I^{n})f\leq c2^{(k-i)n}\int_{\breve{Q}%
_{k,z}}\left\vert f(\nu )\right\vert d\nu ,
\end{equation*}%
where $c$ is independent of $\lambda $ and $\breve{Q}_{k-i,z}=%
\prod_{j=1}^{n}I_{j}$\ with 
\begin{equation*}
I_{j}=\big[2^{i-k}(z_{j}-M-2),2^{i-k}(z_{j}+M+3)\big),\quad j=1,...,n.
\end{equation*}%
Hence%
\begin{equation*}
\bar{\delta}^{M}(x+2^{i-k+1}I^{n})f\lesssim 2^{(k-i)n}\big\|f\big\|%
_{L_{1}(\cup _{l=1}^{T}Q_{k-i,z_{l}})}
\end{equation*}%
for some $T\in \mathbb{N}$. Therefore%
\begin{eqnarray*}
&&\tilde{t}_{k-i,z}\bar{\delta}^{M}(x+2^{i-k+1}I^{n})f \\
&\lesssim &2^{(k-i)n}\sum_{l=1}^{T^{n}}\sum_{h\in \mathbb{Z}%
^{n}:Q_{0,h}\subset Q_{k-i,z_{l}}}\big\|t_{k}(\lambda ^{-1}\cdot )\big\|%
_{L_{p}(Q_{k-i,z})}\big\|f\big\|_{L_{1}(Q_{0,h})}
\end{eqnarray*}%
for any $x\in Q_{k-i,z}$, with $k\in \{1,...,i\}$ and $z\in \mathbb{Z}^{n}$.
Using the fact that $Q_{k-i,z}\subset \cup _{l=1}^{T}Q_{k-i,z_{l}}$ and
Lemma \ref{Ap-Property}/(iii) to estimate this expression by 
\begin{eqnarray*}
&&c2^{(k-i)(n-\frac{n}{\theta })}\sum_{l=1}^{T}\sum_{h\in \mathbb{Z}%
^{n}:Q_{0,h}\subset Q_{k-i,z_{l}}}\bar{t}_{k,h}\big\|f\big\|_{L_{1}(Q_{0,h})}
\\
&\lesssim &2^{(k-i)(n-\frac{n}{\theta })+k\alpha
_{2}}\sum_{l=1}^{T}\sum_{h\in \mathbb{Z}^{n}:Q_{0,h}\subset Q_{k-i,z_{l}}}%
\bar{t}_{0,h}\big\|f\big\|_{L_{1}(Q_{0,h})},
\end{eqnarray*}%
by $\mathrm{\eqref{Asum2}}$ with $\bar{t}_{k,h}=\big\|t_{k}(\lambda
^{-1}\cdot )\big\|_{L_{p}(Q_{0,h})}$, $k\in \{0,1,...,i\}$. Observe that $%
\left( 1+\left\vert x-h\right\vert \right) ^{d}\lesssim 2^{(i-k)d}$ for any $%
x\in Q_{k-i,z}$ and any $d\in \mathbb{N}$, such that $Q_{0,h}\subset
Q_{k-i,z_{l}}$, $l\in \{1,...,T\}$. Therefore%
\begin{equation*}
\tilde{t}_{k-i,z}\bar{\delta}^{M}(x+2^{i-k+1}I^{n})f
\end{equation*}%
can be estimated by%
\begin{equation*}
c2^{(k-i)(n-\frac{n}{\theta }-d)+k\alpha _{2}}\sum_{h\in \mathbb{Z}%
^{n}}\left( 1+\left\vert x-h\right\vert \right) ^{-d}\bar{t}_{0,h}\big\|f%
\big\|_{L_{1}(Q_{0,h})}
\end{equation*}%
for some positive constant $c$ independent of $i,x$ and $k$. Our estimate
use partially some decomposition techniques already used in \cite{FJ90}. For
any $j\in \mathbb{N}$ and any $x\in \mathbb{R}^{n}$, we define 
\begin{equation*}
\Omega _{j}=\{h\in \mathbb{Z}^{n}:2^{j-1}<\left\vert x-h\right\vert \leq
2^{j}\}\text{\quad and\quad }\Omega _{0}=\{h\in \mathbb{Z}^{n}:\left\vert
x-h\right\vert \leq 1\}.
\end{equation*}%
\ Then,%
\begin{eqnarray*}
\sum_{h\in \mathbb{Z}^{n}}\bar{t}_{0,h}\left( 1+\left\vert x-h\right\vert
\right) ^{-d}\big\|f\big\|_{L_{1}(Q_{0,h})} &=&\sum\limits_{j=0}^{\infty
}\sum\limits_{h\in \Omega _{j}}\bar{t}_{0,h}\left( 1+\left\vert
x-h\right\vert \right) ^{-d}\big\|f\big\|_{L_{1}(Q_{0,h})} \\
&\lesssim &\sum\limits_{j=0}^{\infty }2^{-dj}\sum\limits_{h\in \Omega _{j}}%
\bar{t}_{0,h}\big\|f\big\|_{L_{1}(Q_{0,h})}.
\end{eqnarray*}%
For any $0<\kappa <\min (1,p)$, the last expression is bounded by%
\begin{equation}
c\sum\limits_{j=0}^{\infty }2^{-dj}\Big(\sum\limits_{h\in \Omega _{j}}\bar{t}%
_{0,h}^{\kappa }\big\|f\big\|_{L_{1}(Q_{0,h})}^{\kappa }\Big)^{1/\kappa }, 
\notag
\end{equation}%
which can be rewritten as%
\begin{equation}
c\sum\limits_{j=0}^{\infty }2^{(\frac{n}{\kappa }-d)j}\Big(%
2^{-jn}\int\limits_{\cup _{h\in \Omega j}Q_{0,h}}\sum\limits_{m\in \Omega
_{j}}\bar{t}_{0,m}^{\kappa }\big\|f\big\|_{L_{1}(Q_{0,m})}^{\kappa }\chi
_{0,m}(y)dy\Big)^{1/\kappa }.  \label{claim-to prove1}
\end{equation}%
If $y\in \cup _{h\in \Omega _{j}}Q_{0,h}$, then $y\in Q_{0,h}$ for some $%
h\in \Omega _{j}$ and $2^{j-1}<\left\vert x-h\right\vert \leq 2^{j}$. From
this it follows that%
\begin{equation*}
\left\vert y-x\right\vert \leq \left\vert y-h\right\vert +\left\vert
x-h\right\vert \leq \left\vert y-h\right\vert +2^{j}\lesssim \sqrt{n}%
+2^{j}<2^{j+\varrho _{n}},\quad \varrho _{n}\in \mathbb{N},
\end{equation*}%
which implies that $y$ is located in some ball $B(x,2^{j+\varrho _{n}})$. By
taking $d$ large enough, such that $d>\frac{n}{\kappa }$, 
\eqref{claim-to
prove1} does not exceed 
\begin{equation*}
c\mathcal{M}_{\kappa }\Big(\sum\limits_{m\in \mathbb{Z}^{n}}\bar{t}_{0,m}%
\big\|f\big\|_{L_{1}(Q_{0,m})}\chi _{0,m}\Big)(x),\quad x\in \mathbb{R}^{n}.
\end{equation*}%
Clearly and since $0<\kappa <\min (1,p)$, the last\ term in $L_{p}(\mathbb{R}%
^{n})$-quasi-norm is bounded by%
\begin{equation*}
c2^{(k-i)(n-\frac{n}{\theta }-d)+k\alpha _{2}}\Big(\sum_{m\in \mathbb{Z}^{n}}%
\bar{t}_{0,h}^{p}\big\|f\big\|_{L_{1}(Q_{0,m})}^{p}\Big)^{\frac{1}{p}%
}\lesssim 2^{(k-i)(n-\frac{n}{\theta }-d)+k\alpha _{2}}H\big\|f\big\|_{%
\tilde{F}_{p,q}^{M}(\mathbb{R}^{n},\{t_{k}\}_{k\in \mathbb{N}_{0}})}^{\ast }.
\end{equation*}%
Therefore $\mathrm{\eqref{est3}}$, with $\sum\limits_{k=1}^{i}$ in place of $%
\sum\limits_{k=1}^{\infty }$ is bounded by%
\begin{equation*}
C\lambda ^{\alpha _{2}}H\big\|f\big\|_{\tilde{F}_{p,q}^{M}(\mathbb{R}%
^{n},\{t_{k}\}_{k\in \mathbb{N}_{0}})}^{\ast },
\end{equation*}%
since $\alpha _{2}>\max \big(\frac{n}{\theta }-n,\frac{n}{\theta }-\frac{n}{p%
}\big)$ and for some positive constant $C$ independent of $\lambda $. The
proof is complete.\newline

\textbf{Acknowledgements.} The author would like to thank W. Sickel for
valuable suggestions.


\begin{thebibliography}{99}
\bibitem{B03} O.V. Besov, Equivalent normings of spaces of functions of
variable smoothness, Function spaces, approximation, and differential
equations, A collection of papers dedicated to the 70th birthday of Oleg
Vladimorovich Besov, a corresponding member of the Russian Academy of
Sciences, Tr. May. Inst. Steklova, vol. 243, Nauka, Moscow 2003, pp. 87--95;
English transl. in Proc. Steklov Inst. Math. \textbf{243} (2003), 80--88.

\bibitem{B05} O.V. Besov, Interpolation, embedding, and extension of spaces
of functions of variable smoothness, Investigations in the theory of
functions and differential equations, A collection of papers dedicated to
the 100th birthday of academician Sergei Mikhailovich Nikol'skii, Tr. Mat.
Inst. Steklova, vol. 248, Nauka, Moscow 2005, pp. 52--63; English transl.
Proc. Steklov Inst. Math. \textbf{248} (2005), 47--58.

\bibitem{Bui82} H.Q. Bui, Weighted Besov and Triebel spaces: interpolation
by the real method, Hiroshima Math. J. \textbf{12} (1982), 581--605.

\bibitem{BPT97} H.Q. Bui, M. Paluszy\'{n}ski, M.H. Taibleson.
Characterization of the Besov-Lipschitz and Triebel-Lizorkin spaces. The
case $q<1$, J. Fourier Anal. Appl. \textbf{3} (1997), (Spec. Iss.) 837--846.

\bibitem{CF88} F. Cobos, D.L. Fernandez, Hardy-Sobolev spaces and Besov
spaces with a function parameter, In M. Cwikel and J. Peetre, editors,
Function spaces and applications, volume 1302 of Lect. Notes Math., pages
158--170. Proc. US-Swed. Seminar held in Lund, June, 1986, Springer, 1988.

\bibitem{D6} D. Drihem, Besov spaces with general weights, submitted.

\bibitem{D6.1} D. Drihem, Triebel-Lizorkin spaces with general weights,
preprint.

\bibitem{D7} D. Drihem, Spline representations of Lizorkin-Triebel spaces
with general weights, submitted.

\bibitem{ET96} D. Edmunds, H. Triebel, Spectral theory for isotropic fractal
drums. C.R. Acad. Sci. Paris 326, s\'{e}rie I (1998), 1269--1274.

\bibitem{ET99} D. Edmunds, H. Triebel, Eigenfrequencies of isotropic fractal
drums. Oper. Theory: Adv. \& Appl. \textbf{110} (1999), 81--102.

\bibitem{FL06} W. Farkas, H.-G. Leopold, Characterisations of function
spaces of generalised smoothness, Annali di Mat. Pura Appl. \textbf{185}
(2006), 1--62.

\bibitem{FeSt71} C. Fefferman, E.M. Stein, Some maximal inequalities. Amer.
J. Math. \textbf{93} (1971), 107--115.

\bibitem{FJ90} M. Frazier, B. Jawerth,\ A discrete transform and
decomposition of distribution spaces, J. Funct. Anal. \textbf{93} (1990),
34--170.

\bibitem{GR85} J. Garc\.{\i}a-Cuerva, J. L. Rubio de Francia. Weighted Norm
Inequalities and Related Topics, North-Holland Mathematics Studies, 116.
Notas de Matem\'{a}tica [Mathematical Notes], 104. North-Holland Publishing
Co., Amsterdam, 1985.

\bibitem{L. Graf08} L. Grafakos. Classical Fourier Analysis, volume 249 of
Graduate Texts in Mathematics. Springer, New York, 3nd edition, 2014.

\bibitem{Go79} M.L. Goldman, A description of the traces of some function
spaces, Trudy Mat. Inst. Steklov. 150 (1979), 99--127, (Russian) English
transl.: Proc. Steklov Inst. Math. 1981, no. 4 (150).

\bibitem{Go83} M.L. Goldman, A method of coverings for describing general
spaces of Besov type, Trudy Mat. Inst. Steklov. 156 (1980), 47--81,
(Russian) English transl.: Proc. Steklov Inst. Math. 1983, no. 2 156.

\bibitem{Ka83} G.A. Kalyabin, Description of functions from classes of
Besov-Lizorkin-Triebel type, Tr. Mat. Inst. Steklova. 156 (1980), 82--109,
(Russian) English translation: Proc. Steklov Inst. Math. 1983, 156.

\bibitem{Kl87} G.A. Kalyabin, P.I. Lizorkin, Spaces of functions of
generalized smoothness, Math. Nachr. \textbf{133} (1987), 7--32.

\bibitem{Mo01} S.D. Moura, Function spaces of generalised smoothness, Diss.
Math. \textbf{143} (2001).

\bibitem{Mu72} B. Muckenhoupt, Weighted norm inequalities for the Hardy
maximal function, Trans. Amer. Math. Soc. \textbf{165} (1972), 207--226.

\bibitem{HaS08} D.D. Haroske, S.D. Moura, Continuity envelopes and sharp
embeddings in spaces of generalized smoothness, J. Funct. Anal. \textbf{254}%
(6) (2008), 1487--1521.

\bibitem{IzSa12} M. Izuki, Y. Sawano, Atomic decomposition for weighted
Besov and Triebel-Lizorkin spaces, Math. Nachr. \textbf{285} (2012),
103--126.

\bibitem{Naibo} V. Naibo, On the bilinear H\"{o}rmander classes in the
scales of Triebel-Lizorkin and Besov spaces, J. Fourier Anal. Appl. \textbf{5%
} (2015), 1077--1104.

\bibitem{Si} W. Sickel, Spline representations of functions in
Besov-Triebel-Lizorkin spaces on $\mathbb{R}^{n}$, Forum Math. \textbf{2}
(1990), 451--475.

\bibitem{Sc09} C. Schneider, On dilation operators in Besov spaces, Rev.
Mat. Complut. \textbf{22}(1) (2009), 111--128.

\bibitem{ScV09} C. Schneider, J. Vyb\'{\i}ral, On dilation operators in
Triebel-Lizorkin spaces, Funct. Approx. Comment. Math. \textbf{41}(2)
(2009), 139--162.

\bibitem{T1} H. Triebel, Theory of Function Spaces, Birkh\"{a}user Verlag,
Basel, 1983.

\bibitem{T2} H. Triebel, Theory of Function Spaces II, Birkh\"{a}user
Verlag, Basel, 1992.

\bibitem{T3} H. Triebel, Local function spaces, heat and Navier-Stokes
equations. European Math. Soc. Publishing House, Z\"{u}rich, 2013.

\bibitem{Ty14} A.I. Tyulenev, Description of traces of functions in the
Sobolev space with a Muckenhoupt weight, Function spaces and related
questions in analysis, A collection of papers dedicated to the 80th birthday
of Oleg Vladimorovich Besov, a corresponding member of the Russian Academy
of Sciences, Tr. Mat. Inst. Steklova, vol. 284, MAIK, Moscow 2014, pp.
288-303; English transl. in Proc. Steklov Inst. Math. 284 (2014) 280--295.

\bibitem{Ty15} A.I. Tyulenev, Some new function spaces of variable
smoothness, Sbornik Mathematics. \textbf{206} (2015), 849--891.

\bibitem{Ty-N-L} A.I. Tyulenev, Besov-type spaces of variable smoothness on
rough domains, Nonlinear Anal. \textbf{145} (2016), 176--198.

\bibitem{Ty-151} A.I. Tyulenev, On various approaches to Besov-type spaces
of variable smoothness, J. Math. Anal. Appl. \textbf{451} (2017), 371--392.

\bibitem{V08} J. Vyb\'{\i}ral, On dilation operators and sampling numbers,
J. Funct. Spaces Appl. \textbf{6} (2008), 17--46.
\end{thebibliography}
\end{document}